\documentclass[a4paper,10pt]{article}
\usepackage{latexsym} 
\usepackage{amsmath,amssymb}
\usepackage{verbatim}
\usepackage[dvips]{graphicx}

\newtheorem{theorem}{Theorem}[section]
\newtheorem{lemma}{Lemma}[theorem]
\newtheorem{proposition}{Proposition}[theorem]
\newtheorem{corollary}{Corollary}[theorem]
\newtheorem{definition}{Definition}[section]
\newtheorem{preexample}{Example}[section]
\newenvironment{example}{\begin{preexample}}{\end{preexample}}
\newtheorem{preremark}{Remark}
\newenvironment{remark}{\begin{preremark}\rm}{\end{preremark}}

\newcommand{\qed }{ \hfill $\Box$ }

\newcommand{\g}{\mathfrak{g}}
\renewcommand{\log}[1]{ \mathcal{L}^{G\mathfrak{g}}{#1}}
\begin{document}

\begin{center}
{\Large The It\^o exponential on Lie Groups\\}

\end{center}

\vspace{0.3cm}

\begin{center}
{\large  Sim\~ao N. Stelmastchuk}  \\

\textit{Departamento de  Matem\'atica, Universidade Estadual do Paran\'a,\\ 84600-000 -  Uni\~ao da Vit\' oria - PR,
Brazil. e-mail: simnaos@gmail.com}
\end{center}

\vspace{0.3cm}

\begin{abstract}
Let $G$ be a Lie Group with a complete, left invariant connection $\nabla^{G}$. Denote by $\g$ the Lie algebra of $G$ which is equipped with a complete connection $\nabla^{\g}$. Our main goal is to introduce the concept of the It\^o exponential and the It\^o logarithm. As a result, we characterize the martingales in $G$ with respect to the left invariant connection $\nabla^{G}$. Also, assuming that connection function $\alpha:\g \times \g \rightarrow \g$ associated to $\nabla^{G}$ satisfies $\alpha(M,M) = 0$ for all $M \in \g$ we obtain a stochastic Campbell-Hausdorff formula. Further, from this stochastic Campbell-Hausdorf formula we present a way to construct martingales in Lie group. In consequence, we show that a product of harmonic maps with value in $G$ is a harmonic map. To end, we apply this study in some matrix Lie groups.
\end{abstract}


\vspace{0.3cm} \noindent {\bf MSC2010 subject classification:}  22E99, 53C43, 58E20, 58J65, 60H30

\section{Introduction}

M. Hakim-Dowek and D. L\'epingle introduced the exponential stochastic in Lie Groups by first time in \cite{hakim}.  Their idea can be interpreted in the following way. Let $G$ be a Lie group and $\g$ its Lie algebra. Given a semimartingale $M$ in $\g$, the stochastic exponential  $X=e(M)$ is the solution, in the Stratonovich sense, of the stochastic differential equation
\[
\delta X = L_{X*} \delta M, \\ X_{0}=e.
\]
Furthermore, they showed a existence of the inverse of the stochastic exponential, which is known as the stochastic logarithm. In this paper \cite{hakim} was developed a serial of result about the stochastic exponential and the stochastic logarithm.

The concept of the stochastic exponential and the stochastic logarithm in Lie group has been studied and applied in some situations. For example, M. Arnaudon developed studies of  the stochastic exponential in Lie groups in the case that $G$ has a left invariant connection \cite{arnaudon1}. He also used the stochastic exponential to study the martingales and Brownian motions in homogeneous space \cite{arnaudon2}. The characterization of  the semimartingales, martingales and Brownian motions in a principal fiber bundle, due to M. Arnaudon and S. Paycha \cite{arnaudon3}, is obtained with the stochastic exponential. In the case that a Lie Group $G$ has a bi-invariant metric, P. Catuogno and P. Ruffino in \cite{catuogno1} used the stochastic exponential and the stochastic logarithm to show that the product of harmonic maps with values in $G$ is a harmonic map.

The developed of the stochastic exponential has been done without to takes in account some geometry of the Lie groups. Despite of  the above studies have worked with some types of connections. In a nutshell, this fact occurs because the integral of Stratonovich does not have intrinsically the geometry of the smooth manifolds.  Unlike, the integral of It\^o on a smooth manifold intrinsically has the information of the geometry of the smooth manifold. Considering this fact  we introduce a stochastic exponential and a stochastic logarithm in the It\^o sense.

Let $\nabla^{G}$ be a complete, left invariant connection on $G$ and $\nabla^{\g}$ a complete connection on $\g$. The It\^o exponential and the It\^o logarithm are the solutions of the following stochastic differential equations, respectively,
\[
d^{G}X = L_{X*}d^{\g}M, \ \ X_{0}=e
\]
and
\[
d^{g}N = L_{Y^{-1}*}d^{G}Y, \ \ Y_{0} = 0.
\]
We denote this solution by $e^{G\g}(M)$ and $\log{(Y)}$, respectively.

Our first work is to show that this equations have unique solutions that do not explode in a finite time, since $\nabla^{G}$ and $\nabla^{\g}$ are completes. Also, we show that the operators $e^{G\g}$ and $\log{}$ are inverses. As a result, we get that every $\nabla^{G}$-martingales is given by $e^{G\g}(M)$ for a $\nabla^{\g}$-martingale $M$ in $\g$.

An observable fact is that the It\^o exponential and It\^o logarithm are dependents of the connections $\nabla^{G}$ and $\nabla^{\g}$, however it is not  necessary a kind of the correspondence between these connections.

Other work is to construct a stochastic Campbell-Hausdorf formula. Given a left invariant connection $\nabla^{G}$ there exists a unique bilinear form $\alpha$ on $\g$ associated to $\nabla^{G}$ (see for example \cite{helgason}). If we suppose that $\alpha(M,M) = 0$ for all $M \in \g$, then $\nabla^{G}$ is complete. Under this hypothesis, with property of null quadratic variation property, which is defined in \ref{campbellnullprop1}, we show a stochastic Campbell-Hausdorff formula. This formula help us to create a way to construct martingales in the Lie group with respect to $\nabla^{G}$. A little bit, from this result we can show that a product of harmonic maps with values in $G$ is a harmonic map.

To end, we show a relation between It\^o logarithm and stochastic logarithm. This relation will give us the geometry of Lie group in therms of the It\^o logarithm. We apply this result in the study of the martingales in some matrix Lie group equipped with a class of the left invariant metric, namely, the Euclidian motion group $SE(3)$ and the three-dimensional non compact Lie Groups $SE(2)$, $E(1,1)$, $N^{3}$ and $SL(2,\mathbb{R})$.

%
%
\section{Preliminaries}

In this work we use freely the concepts and notations of  P. Protter \cite{Protter}, E. Hsu \cite{hsu}, P. Meyer
\cite{Mey3}, M. Emery \cite{emery1} and \cite{emery3}, and S. Kobayashi and N. Nomizu \cite{kobay}. We suggest the
reading of \cite{catuogno2} for a complete survey about the objects of this section. From now on the adjective smooth
means $C^{\infty}$.

Let $(\Omega, \mathcal{F},(\mathcal{F}_{t})_{t\geq0}, \mathbb{P})$ be a probability space which satisfies the usual hypotheses (see for example \cite{emery1}). Our basic assumption is that every stochastic process is continuous.

Let $M$ be a smooth manifold and $X_{t}$ a continuous stochastic process with values in $M$. We call $X_{t}$ a semimartingale if, for all $f$ smooth function, $f(X_{t})$ is a real semimartingale.

Let $M$ be a smooth manifold endowed with a symmetric connection $\nabla^{M}$. From now on we make the assumption: {\it all connections will be symmetrics}. Let $X$ be a semimartingale in $M$ and $\theta$ a 1-form on $M$ defined along $X$. Let $(U, x_{1},\ldots, x_{n})$ be a local coordinate system on $M$. We define the Stratonovich and the It\^o integrals, respectively, of $\theta$ along $X$, locally, by
\[
\int_0^t\theta \delta X_{s} = \int_{0}^{t} \theta_{i}(X_{s})dX^{i}_{s} + \frac{1}{2}\int_{0}^{t}
\frac{\partial\theta_{i}}{\partial x^{j}}(X_{s})d[X^{i},X^{j}]_{s},
\]
and
\begin{equation}\label{itointegral}
\int_0^t\theta d^{M} X_{s} = \int_{0}^{t} \theta_{i}(X_{s})dX^{i}_{s} + \frac{1}{2}\int_{0}^{t}
\Gamma_{jk}^{i}(X_{s})\theta_{i}(X_{s})d[X^{j},X^{k}]_{s},
\end{equation}
where $\theta = \theta_{i} dx^{i}$ with $\theta_{i}$ smooth functions and $\Gamma_{jk}^{i}$ are the Christoffel symbols
of the connection $\nabla^{M}$. Let $b \in T^{(2,0)}M$ be defined along $X$. We define the quadratic integral on $M$
along $X$, locally, by
\[
\int_0^{t}b\;(dX,dX)_{s} = \int_{0}^{t} b_{ij}(X_{s})d[X^{i},X^{j}]_{s},
\]
where $b = b_{ij} dx^{i} \otimes dx^{j}$ with  $b_{ij}$ smooth functions.

A direct consequence of the definitions above is the Stratonovich-It\^o formula of conversion given by
\begin{equation}\label{conversion}
\int_0^t\theta \delta X_{s} = \int_0^t\theta d^{M} X_{s}+\frac{1}{2}\int_0^t\nabla^{M}\theta\;(dX,dX)_{s}.
\end{equation}

A semimartingale $X$ with values in $M$ is called a $\nabla^{M}$-martingale if $\int \theta d^{M} X$ is a real local martingale for all $\theta \in \Gamma(TM^*)$.

Let $M$ be a Riemannian manifold with a metric $g$. A semimartingale $B$ in $M$ is said a $g$-Brownian motion if $B$ is a $\nabla^{g}$-martingale, being $\nabla^{g}$ the Levi-Civita connection of $g$, and for any section $b$ of $T^{(2,0)}M$ we have
\begin{equation}\label{Brownian}
\int_0^tb(dB,dB)_{s}=\int_0^t \mathrm{tr}\,b(B_s)ds.
\end{equation}

Let $M$ and $N$ be smooth manifolds with connections $\nabla^{M}$ and $\nabla^{N}$, respectively, $\theta$ a section of $TN^*$, $b$ a section of $T^{(2,0)}N$ and $F:M\rightarrow N$ a smooth map. The geometric It\^o formula is given by
\begin{equation}\label{itoformula}
 \int_{0}^{t} \theta d^{N}F(X_{s}) = \int_{0}^{t} F^{*}\theta d^{M}X_{s} + \frac{1}{2} \int_{0}^{t} \beta_{F}^{*}\theta(dX,dX)_{s},
\end{equation}
where $\beta_{F}$ is the second fundamental form of $F$ with respect to $\nabla^{M}$ and $\nabla^{N}$.

%
%
\section{The It\^o exponential and It\^o logarithm}

Let $G$ be a Lie Group and $\mathfrak{g}$ its Lie algebra. Let us denote by $L_{g}$ the left translation on $G$. From this we can construct the following family of linear applications on $\mathfrak{g}^{*} \otimes TG$: since $\mathfrak{g}$ is isomorphic to $T\mathfrak{g}$, we consider that the left translation is a linear application $L_{g*}(e):\mathfrak{g} \rightarrow TG$ for every $g \in G$. We observe that this family of applications is smooth in the following sense. Taking $E \in \mathfrak{g}$ we obtain a smooth left invariant vector field $X \in TG$ such that $L_{g*}(e)(E) = X_{g}$. Therefore $L_{g*}(e)$ is a smooth family from $\mathfrak{g} \times G$ into $TG$ (see for instant Definition 6.34 in \cite{emery1}).

We endow $G$ with a left-invariant connection $\nabla^{G}$ and $\mathfrak{g}$ with a connection $\nabla^{\mathfrak{g}}$.Let $X$ be a semimartingale in $G$ and $M$ a semimartingale in $\mathfrak{g}$. One says that $X$ is a solution to the stochastic differential equation
\begin{equation}\label{exponential}
d^{G}X_{t} = L_{X_{t}*}(e) d^{\g}M_{t},
\end{equation}
if, for  every 1-form $\theta$ on $G$, the real semimartingales $\int_{0}^{t} \theta d^{G}X_{s}$ and \linebreak
$\int_{0}^{t} L_{X_{s}*}(e) \theta d^{\g}M_{s}$ are equal.

Firstly, one may observe that the solution of the stochastic differential equation above is invariant because of the left invariance of the connection $\nabla^{G}$.

\begin{proposition}\label{ExpLogGroupprop1}
Let $G$ be a Lie group and $\mathfrak{g}$ its Lie algebra. Assume that $\nabla^{G}$ is a left-invariant connection on $G$ and $\nabla^{\mathfrak{g}}$ is a connection on $\mathfrak{g}$. Suppose that $Y_{t}$ is a solution of (\ref{exponential}). If $\xi$ is a random variable with values in $G$, then $X_{t}=\xi Y_{t}$ is also a solution of  (\ref{exponential}).
\end{proposition}
\begin{proof}
We begin denoting the product on Lie group $G$ by $m$. Let $\theta$ be a 1-form on $G$. As a function to $m$, the integral of It\^o along $X_{t}$ is writing by
\[
\int_{0}^{t} \theta d^{G}X_{s} = \int_{0}^{t} \theta d^{G}\xi Y_{s} = \int_{0}^{t} \theta d^{G}m(\xi,Y_{s}).
\]
The geometric It\^o formula (\ref{itoformula}) get
\[
\int_{0}^{t} \theta d^{G}X_{s} = \int_{0}^{t} m^{*}\theta d^{(G\times G)}(\xi ,Y_{s}) + \frac{1}{2}\int_{0}^{t}
\beta_{m}^{*}\theta(d(\xi,Y_{s}),d(\xi,Y_{s})).
\]
From Proposition 3.15 in \cite{emery3} we see that
\[
\int_{0}^{t} \theta d^{G}X_{s} = \int_{0}^{t} (R_{Y_{s}}^{*}\theta) d^{G}\xi + \int_{0}^{t} (L_{\xi}^{*}\theta) d^{G}Y_{s} + \frac{1}{2}\int_{0}^{t} \beta_{m}^{*}\theta(d(\xi,Y_{s}),d(\xi,Y_{s})).
\]
We see that $\xi$ is a constant process, consequently,
\[
\int_{0}^{t} \theta d^{G}X_{s} = \int_{0}^{t} (L_{\xi}^{*}\theta) d^{G}Y_{s} + \frac{1}{2}\int_{0}^{t} \beta_{m}^{*}\theta(d(\xi,Y_{s}),d(\xi,Y_{s})).
\]
We claim that the $\beta_{m}(d(\xi,Y_{t}),d(\xi,Y_{t}))$ is null. In fact, take $0 \in T_{g}G$ and a left invariant vector field $Y$ on $G$. Here, $0$ is the vector associated to the constant process $\xi$. Then
\begin{eqnarray*}
\beta_{m}(0,Y)
& = & \nabla^{G}_{m_{*}(0,Y)}m_{*}(0,Y) - m_{*}\nabla^{G \times G}(0, Y)\\
& = & \nabla^{G}_{R_{h*} 0 + L_{g*}(Y)}(R_{h*} 0 + L_{g*}(Y)) - m_{*}\nabla^{G \times G}(0, Y)\\
& = & \nabla^{G}_{L_{g*}Y}L_{g*}Y - L_{g*}(\nabla^{G}_{Y} Y)\\
& = & L_{g*}(\nabla^{G}_{Y}Y) - L_{g*}(\nabla^{G}_{Y} Y)\\
& = & 0,
\end{eqnarray*}
where in forth equality we use the fact that $\nabla^{G}$ is a left invariant connection. Thus we get
\[
\int_{0}^{t} \theta d^{G}X_{s} = \int_{0}^{t} (L_{\xi}^{*}\theta) d^{G}Y_{s}.
\]
As $Y_{t}$ is a solution of (\ref{exponential}) we have
\[
 \int_{0}^{t} \theta d^{G}X_{s} = \int_{0}^{t} L_{Y_{s}}^{*}(e)L_{\xi}^{*}(Y_{s})\theta d^{\g}M_{s}.
\]
This gives
\[
 \int_{0}^{t} \theta d^{G}X_{s} = \int_{0}^{t} L_{X_{t}}^{*}(e)\theta d^{\g}M_{s}.
\]
Therefore we conclude that $X_{t}$ is a solution of (\ref{exponential}).\qed
\end{proof}

The idea to show the existence and unicity of the solution of (\ref{exponential}) is to construct a second order stochastic differential equations from this and use the results of the existence and unicity given by \cite{emery1}.

\begin{proposition}\label{ExpLogGroupteo1}
Let $M$ be a semimartingale in $\mathfrak{g}$ and $X_{0}$ a $\mathcal{F}_{0}$-measurable random variable in $G$. There exist a predictable stopping time $\zeta$ and a $G$-valued semimartingale $X$ in $G$ on the interval $[ 0, \zeta [$, with initial condition $X_{0}$, solution to $(\ref{exponential})$ and exploding to times $\zeta$ on the event $\{\zeta < \infty \}$. Moreover, the following uniqueness and maximality properties holds: if $\zeta'$ is a predictable time and
$X'$ a solution starting from $X_{0}$ defined on $[0, \zeta'[$, then $\zeta' \leq \zeta$ and $X' = X$ on $[0, \zeta'[$.
\end{proposition}
\begin{proof}
From the It\^o transfer principle, see Theorem 12 in \cite{emery2}, the stochastic differential equation (\ref{exponential}) is equivalent to the intrinsic second order differential equation $d^{2}X = f(X,M) d^{2}M$, where $f:\tau_{M_{t}(\omega)}\g \rightarrow \tau_{X_{t}(\omega)}G$ is the unique semi-affine Schwartz morphism with  $L_{X_{t}(\omega)*}(e)$ as restriction to the first order. From Lemma 11 in \cite{emery2} we see that $f$ is a family of Schwartz morphism which depend smoothly upon $(e,g)$, for all $g \in G$. From Theorem 6.41 in \cite{emery1}, with a $\mathcal{F}_{0}$-measurable random variables $X_{0}$ in $G$ as an initial condition, there exists a predictable stopping time $\zeta >0$ and a $G$-valued semimartingale $X$ on the interval $[ 0, \zeta [ $ which  is the solution of  $d^{2}X = f(X,M) d^{2}M$ and exploding at time $\zeta$ on the event $\{\zeta < \infty\}$. Moreover,
the uniqueness and maximality properties holds in the following sense: if $\zeta'$ is a predictable time and $X'$ is a solution starting from $X_{0}$ defined on $[0, \zeta'[$, then $\zeta' \leq \zeta$ and $X' = X$ on $[0, \zeta'[$. Again, It\^o
transfer principle assure the unique and maximal solution of (\ref{exponential}) with initial condition $X_{0}$ and predictable stopping time $\zeta$.\qed
\end{proof}

This Proposition deals with hypothesis that $\nabla^{G}$ and $\nabla^{\g}$ are any connections on G and g, respectively. However, if we wish study the martingales in $G$, we need some restriction in these connections. This restriction is that connections are completes as one see in the next proposition.

\begin{proposition}\label{ExpLogGroupprop2}
If $X_{t}$ is the solution of the stochastic differential equation (\ref{exponential}), then this life time is infinity, since $\nabla^{G}$ and $\nabla^{\g}$ are completes.
\end{proposition}
\begin{proof}
It is a direct consequence of Theorem 14 in \cite{emery2} and the fact that every geodesic in $G$ and $\g$ with respect to $\nabla^{G}$ and $\nabla^{\g}$, respectively, are extended to all time. \qed
\end{proof}

Proposition \ref{ExpLogGroupprop1} says that we can consider the solution of stochastic differential equation (\ref{exponential}) with initial value $X_{0} =e$ rather than any random variable on $G$. In other side, Proposition \ref{ExpLogGroupprop2} shows that solution of (\ref{exponential}) is in interval $[0,\infty[$. From these facts we give the following definition.

\begin{definition}
Suppose that $\nabla^{G}$ is a complete, left invariant connection on $G$ and $\nabla^{\g}$ is a complete connection on $\g$. We will denote by $e^{G \mathfrak{g}}(M)$ the solution of (\ref{exponential}) with initial condition $X_{0}=e$ and we call it \emph{It\^o stochastic exponential with respect to $\nabla^{G}$ and $\nabla^{\mathfrak{g}}$}.
\end{definition}

In the follow, for simplicity, we will call $e^{G \mathfrak{g}}(M)$ by It\^o exponential.

\begin{remark}
It is well known that the left-invariant connections on $G$ are in one-one correspondence with bilinear forms on $\mathfrak{g}$, see Proposition 1, chapter 3 in \cite{helgason}. However the stochastic differential equation (\ref{exponential}) do not preserve this fact, that is, it is not necessary that $\nabla^{G}$ and  $\nabla^{g}$ have some association.
\end{remark}

As an immediate consequence from Proposition \ref{ExpLogGroupteo1} we have a characterization of the $\nabla^{G}$-martingales in the Lie group $G$.

\begin{corollary}\label{expmart}
The It\^o exponential $e^{G\mathfrak{g}}(M)$ is a $\nabla^{G}$-martingale in $G$ if and only if $M$ is a $\nabla^{\mathfrak{g}}$-martingale in $\mathfrak{g}$.
\end{corollary}

The It\^o exponential yields a semimartingale in $G$ from a semimartingale in $\mathfrak{g}$. M. Hakim-Dowek  and D. L\'epingle \cite{hakim} define, in Stratonovich sense, an inverse of the stochastic exponential, which they called the stochastic logarithm. We wish to get an analogous in the It\^o sense. For this we consider a left-invariant connection $\nabla^{G}$ on $G$ and a connection $\nabla^{\mathfrak{g}}$ on $\mathfrak{g}$. Our idea is to create the process inverse of the It\^o exponential as solution of the following stochastic differential equation
\begin{equation}\label{logaritmo}
d ^{\g}M_{t} = L_{(X_{t})^{-1}*}(X_{t}) d^{G}X_{t}.
\end{equation}
The solution to this stochastic differential equation means that for every 1-form $\psi$ on $\mathfrak{g}$ the real semimartingales $\int \psi d^{\g}M$ and $\int L_{X_{t}^{-1}*}(X_{t}) \psi d^{G}X$ are equal.

An important invariance property to the solution of stochastic differential equation above is obtained if we ask the left invariance property to $\nabla^{G}$.

\begin{proposition}\label{ExpLogGroupprop3}
Let $G$ be a Lie group and $\mathfrak{g}$ its Lie algebra. Assume that $\nabla^{G}$ is a left-invariant connection on $G$ and $\nabla^{\mathfrak{g}}$ is a connection on $\mathfrak{g}$. Suppose that $M_{t}$ is a solution of (\ref{logaritmo}) with respect to a semimartingale $Y_{t}$. If $\xi$ is a random variable in $G$, then $M_{t}$ is a solution of (\ref{logaritmo}) with respect to $X_{t}=\xi Y_{t}$.
\end{proposition}
\begin{proof}
Let $\psi$ be a 1-form on $\mathfrak{g}$. By definition, the It\^o stochastic differential equation (\ref{logaritmo}) means that
\[
\int_{0}^{t} \psi d^{\g}M_{s} = \int_{0}^{t} L_{Y_{s}^{-1}*}\psi d^{G}Y_{s}.
\]
Since $Y_{t} = L_{\xi^{-1}}X_{t}$, it follows that
\[
\int_{0}^{t} \psi d^{\g}M_{s} = \int_{0}^{t} L_{(\xi^{-1}X_{s})^{-1}*}\psi d^{G}L_{\xi^{-1}}X_{s}.
\]
As in the proof of Proposition \ref{ExpLogGroupprop1} we have
\[
\int_{0}^{t} \psi d^{\g}M_{s} = \int_{0}^{t} L_{\xi^{-1}}^{*}L_{(\xi^{-1}X_{s})^{-1}}^{*}\psi d^{G}X_{s}.
\]
We thus get
\[
\int_{0}^{t} \psi d^{\nabla^{\mathfrak{g}}}M_{s} = \int_{0}^{t} L_{X_{s}^{-1}}^{*}\psi d^{\nabla^{G}}X_{s}.
\]
By definition, $M_{t}$ is also a solution of (\ref{logaritmo}) with respect to $X_{t}$. \qed
\end{proof}

The solution and uniqueness to (\ref{logaritmo}) are assured to follow.

\begin{proposition}\label{ExpLogGroupteo2}
Let $X$ be a semimartingale in $G$ and $M_{0}$ a $\mathcal{F}_{0}$-measurable random variable with value in $\mathfrak{g}$. There exists a predictable stopping time $\eta$ and a $\mathfrak{g}$-valued semimartingale $M$ in $\mathfrak{g}$ on the interval $[ 0, \eta [$, with initial condition $M_{0}$, solution to (\ref{logaritmo}) and exploding to times $\eta$ on the event $\{\eta < \infty \}$. Moreover, the following uniqueness and maximality properties hold: if $\eta'$ is a predictable time and $M'$ a solution starting from $M_{0}$ defined on $[0, \eta'[$, then $\eta' \leq \eta$ and $M' = M$ on $[0, \eta'[$.
\end{proposition}
\begin{proof}
The proof is analogous to one in Proposition \ref{ExpLogGroupteo1}. \qed
\end{proof}

\begin{proposition}\label{ExpLogGroupprop4}
The solution of (\ref{logaritmo}) has a life time infinity, since $\nabla^{G}$ and $\nabla^{g}$ are completes.
\end{proposition}
\begin{proof}
The same argument that is used to prove Proposition \ref{ExpLogGroupprop2}. \qed
\end{proof}

Propositions \ref{ExpLogGroupprop3} - \ref{ExpLogGroupprop4} yield the well definition of the It\^o logarithm.

\begin{definition}
Suppose that $\nabla^{G}$ is a complete, left invariant connection and $\nabla^{\g}$ is a complete connection. The solution to (\ref{logaritmo}), with initial condition $M_{0}=0$, is called \emph{It\^o stochastic logarithm with respect to $\nabla^{G}$ and $\nabla^{\mathfrak{g}}$} and it is denoted by $\log{(X)}$.
\end{definition}

For simplicity, we call $\log{(X)}$ by It\^o logarithm. Now, one may see, as in the case of the It\^o exponential, the clear relation between $\nabla^{\mathfrak{g}}$-martingales and $\nabla^{G}$-martingales.

\begin{corollary}\label{ExpLogGroupcor2}
A semimartingale $X_{t}$ in $G$ is a $\nabla^{G}$-martingale if and only if $\log{(X)}$ is a $\nabla^{\mathfrak{g}}$-martingale in $\g$.
\end{corollary}

Our main intention to work with the It\^o Logarithm  is that it is the inverse of the It\^o exponential.

\begin{theorem}\label{ExpLogGroupteo3}
Let $\nabla^{G}$ be a complete, left invariant connection on $G$ and $\nabla^{\g}$ a
complete connection on $\g$. Let $X_{t}, M_{t}$ be semimartingales in $G$ and $\mathfrak{g}$, respectively. Then
\[
\mathcal{L}^{G\mathfrak{g}}(e^{G\mathfrak{g}}(M_{t})) = M_{t}
\]
and
\[
e^{G\mathfrak{g}}(\mathcal{L}^{G\mathfrak{g}}(X_{t})) = X_{t}.
\]
\end{theorem}
\begin{proof}
Let $\psi$ be a 1-form on $\mathfrak{g}$ and $M_{t}$ a semimartingale in $\mathfrak{g}$. The It\^o exponential is the
semimartingale $e^{G \mathfrak{g}}(M_{t})$ in  $G$. By definition, the It\^o logarithm apply to $e^{G\mathfrak{g}}(M_{t})$ means that
\[
\int_{0}^{t} \psi d^{\g} \mathcal{L}^{G\mathfrak{g}}(e^{G\mathfrak{g}}(M_{s})) = \int_{0}^{t}
L_{e^{G\mathfrak{g}}(M_{s})^{-1*}}\psi d^{G}e^{G\mathfrak{g}}(M_{s}).
\]
Since $e^{G \mathfrak{g}}(M_{t})$ is solution of (\ref{exponential}), it follows that
\[
\int_{0}^{t} \psi d^{\g} \mathcal{L}^{G \mathfrak{g}}(e^{G \mathfrak{g}}(M_{s}))
= \int_{0}^{t} L_{e^{G\mathfrak{g}}(M_{s})^{*}}L_{e^{G\mathfrak{g}}(M_{s})^{-1*}}\psi d^{\g}M_{s}
=  \int_{0}^{t} \psi d^{\g}M_{s}.
\]
As $\psi$ is an arbitrary 1-form on $\mathfrak{g}$ we have $\mathcal{L}^{G\mathfrak{g}}(e^{G\mathfrak{g}}(M_{t})) = M_{t}$. Similarly, for a semimartingale $X_{t}$ in $G$, we get $e^{G\mathfrak{g}}(\mathcal{L}^{G\mathfrak{g}}(X_{t})) = X_{t}$ .\qed
\end{proof}

To end this section, we characterize the martingales in Lie groups.

\begin{theorem}\label{ExpLogGroupteo4}
Let $\nabla^{G}$ be a complete, left invariant connection on $G$ and $\nabla^{\g}$ a complete connection on $\g$. Every $\nabla^{G}$-martingale in $G$ is write as $e^{G\mathfrak{g}}(M_{t})$ for a \linebreak $\nabla^{\mathfrak{g}}$-martingale $M_{t}$ in $\mathfrak{g}$.
\end{theorem}
\begin{proof}
Let $X_{t}$ be a $\nabla^{G}$-martingale in $G$. From Corollary \ref{ExpLogGroupcor2} we see that $\mathcal{L}^{G\mathfrak{g}}(X_{t})$ is a $\nabla^{\mathfrak{g}}$-martingale in $\g$. Taking $M_{t}=\mathcal{L}^{G\mathfrak{g}}(X_{t})$ we conclude that \linebreak $X_{t} = e^{G\mathfrak{g}}(M_{t})$, which
follows from Theorem \ref{ExpLogGroupteo3}. \qed
\end{proof}

\section{Campbell-Hausdorff formulas}

In this section, we give a stochastic Campbell-Hausdorff formula. For this end, we ask for  a condition over the left invariant connection $\nabla^{G}$. More exactly, to the connection function $\alpha: \g \times \g \rightarrow \g$ associated to $\nabla^{G}$ we suppose that $\alpha(A,A) = 0$ for all $A \in \g$. It is known that $\nabla^{G}$ is complete. Before we show the stochastic Campbell-Hausdorff formulas we need to introduce the null quadratic variation property.

\begin{definition}
Let $N$ be a smooth manifold and $X_t,Y_t$ be two semimartingales in $N$. We say that $X_t$ and $Y_t$ have the null quadratic variation property if for any local coordinates system $(x^{1}, \ldots , x^{n})$ on $N$ we have $[X^{i}_t, Y^{j}_t]=0$, for $i,j=1,\ldots, n$, where $X^{i}_t=x^{i} \circ X_t$ and $Y^{j}_t = y^{j} \circ Y_t$.
\end{definition}

\begin{example}
Any two independents semimartingales in a smooth manifold $N$ have the null quadratic variation property.
\end{example}

The null quadratic variation property is good because it is held by the It\^o logarithm.

\begin{proposition}\label{campbellnullprop1}
Let $G$ be a Lie group with a complete, left invariant connection $\nabla^{G}$ and $\g$ its Lie algebra endowed with a complete connection $\nabla^{\g}$. Given two semimartingales  $X_t,Y_t$ in $G$, then $X_t,Y_t$ have the null quadratic variation property if and only if $\log{(X_t)}$ and $\log{(Y_t)}$ have the null quadratic variation property.
\end{proposition}
\begin{proof}
Suppose that $X_t,Y_t$ are two semimartingales in $G$ such that $X_t,Y_t$ have the null quadratic variation property. Thus, for any local coordinates system $(U, x^{1}, \ldots , x^{n})$ on $G$ we see that $[X^{i}_t, Y^{j}_t]=0$, where $X^{i}_t=x^{i}\circ X_t$ and $Y^{j}_t=x^{j}\circ Y_t$, $i,j=1, \ldots, n$. It is sufficient to prove that $[\log{(X_t)}^{\alpha},\log{(Y_t)}^{\beta}] = 0$ for a global coordinate system $(y^{1},\ldots , y^{n})$ on $\g$. By Proposition 7.8 in \cite{emery1},
\[
[\log(X_{t})^{\alpha}, \log(Y_{t})^{\beta}] = [\int_{0}^{t} dy^{\alpha} d\log{X_{s}},\int_{0}^{t} dy^{\beta} d\log{Y_{s}}].
\]
From definition of the It\^o logarithm we see that
\[
[\log(X_{t})^{\alpha}, \log(Y_{t})^{\beta}] = [\int_{0}^{t} dy^{\alpha} L_{X_{s}*}^{-1}(X_{s}) d^{G}X_{s},\int_{0}^{t} dy^{\beta} L_{Y_{s}*}^{-1}(Y_{s}) d^{G}Y_{s}].
\]
Applying the definition of the  integral of It\^o (\ref{itointegral}) yields
\[
[\log(X_{t})^{\alpha}, \log(Y_{t})^{\beta}] = [\sum_{l=1}^{n} \int_{0}^{t}\! (dy^{\alpha} L_{X_{s}*}^{-1}(X_{s}))^{l} dX^{l}_{s}, \sum_{l=1}^{n} \int_{0}^{t}\! (dy^{\beta} L_{X_{s}*}^{-1}(X_{s}))^{k} dX^{k}_{s}].
\]
Interchanging the integral of It\^o with quadratic variation we obtain
\[
[\log(X_{t})^{\alpha}, \log(Y_{t})^{\beta}] = \sum_{l,k=1}^{n}\int_{0}^{t} (dy^{\alpha} L_{X_{s}*}^{-1}(X_{s}))^{l}(dy^{\beta} L_{Y_{s}*}^{-1}(Y))^{k}d[X^{l}_{s},Y^{k}_{s}].
\]
Since $X,Y$ are null quadratic variation, $[\log(X_{t})^{\alpha}, \log(Y_{t})^{\beta}]=0$. It gives the null
quadratic variation property for $\log{(X)}$ and $\log{(Y)}$.

Similarly, one can show that if $\log{(X)}$ and $\log{(Y)}$ have the null quadratic variation property, then $X,Y$ also have the one property. \qed
\end{proof}

The Campbell-Hausdorff formula is then given by

\begin{theorem}\label{campbellnullteo1}
Let $G$ be a Lie group with a complete, left invariant connection $\nabla^{G}$ such that its connection function $\alpha$ satisfies $\alpha(A,A) = 0$ for all $A \in \g$ and $\g$ its Lie algebra endowed with a complete connection $\nabla^{\g}$. Given two semimartingales $M,N$ in $\g$ which satisfy the null quadratic variation property, then
\begin{equation}\label{campbellnullteo1eq1}
e^{G \g}(M_{t} + N_{t}) = e^{G \g}\left(\int_{0}^{t} Ad(e^{G\g}(N_{s}))dM_{s}\right)e^{G\g}(N_{t}).
\end{equation}
For two semimartingales $X,Y$ in $G$, which have the null quadratic variation property, we have
\begin{equation}\label{campbellnullteo1eq2}
\mathcal{L}^{G\g}(X_{t}Y_{t}) = \int_{0}^{t} Ad(Y^{-1}_{s})d\mathcal{L}^{G\g}(X_{s})+ \mathcal{L}^{G\g}(Y_{t}).
\end{equation}
\end{theorem}
\begin{proof}
We begin introducing the following notation
\begin{equation}\label{campbellnullteo1eq3}
 X_{t} = e^{G \g}\left(\int_{0}^{t} Ad(e^{G\g}(N_{s})) dM_{s}\right) \ \ \textrm{and} \ \ \ Y_{t} = e^{G\g}(N_{t}) .
\end{equation}
The proof of (\ref{campbellnullteo1eq1}) is complete if for each left invariant 1-form $\theta$ on $G$
\[
\int_{0}^{t} \theta d^{G}(X_{s}Y_{s}) = \int_{0}^{t} \theta L_{(X_{s}Y_{s})*}(e)d(M_{s}+N_{s}).
\]
Consider the product on Lie group as the application $m: G\times G \rightarrow G$. Using the geometric It\^o formula
(\ref{itoformula}) we get
\begin{eqnarray*}
\int_{0}^{t} \theta d^{G}(X_{s}Y_{s}) & = & \int_{0}^{t} \theta d^{G}m(X_s,Y_s)\\
 & = & \int_{0}^{t} m^{*}\theta d^{G \times G}(X_{s},Y_s) + \frac{1}{2}\int_{0}^{t} \beta_{m}^{*}(d(X_s,Y_s),d(X_s,Y_s)).
\end{eqnarray*}
We have  $\frac{1}{2}\int_{0}^{t} \beta_{m}^{*}(d(X_s,Y_s),d(X_s,Y_s))= 0$, because $\alpha(A,A) = 0$ for all $A \in \g$ and $X,Y$ have the null quadratic variation property. Hence
\[
\int_{0}^{t} \theta d^{G}(X_tY_t) = \int_{0}^{t} m^{*}\theta d^{G \times G}(X_{s},Y_{s}).
\]
From Proposition 3.7 in \cite{emery3} it may be conclude that
\[
\int_{0}^{t} \theta d^{G}(X_{s}Y_{s}) = \int_{0}^{t} R_{Y_{s}}^{*}\theta d^GX_{s} + \int_{0}^{t} L_{X_{s}}^{*}\theta d^GY_{s}.
\]
Replacing (\ref{campbellnullteo1eq3}) in this equality yields
\[
\int_{0}^{t} \theta d^{G}(X_{s}Y_s) = \int_0^t R_{Y_s}^{*}\theta L_{X_s*}(e)Ad(Y_s)dM_s + \int_0^t L_{X_s}^{*}\theta L_{Y_s*}(e) dN_s.
\]
Here, an easy computation shows that
\begin{eqnarray*}
\int_{0}^t \theta d^{G}(X_s Y_s)
& = &  \int_0^t \theta L_{X_sY_s*}(e)dM_s + \int \theta L_{X_sY_s*}(e)dN_s\\
& = & \int_0^t \theta L_{X_sY_s*}(e)d(M_s + N_s),
\end{eqnarray*}
and the proof is complete.

The equality (\ref{campbellnullteo1eq2}) is a direct consequence of (\ref{campbellnullteo1eq1}). \qed
\end{proof}

\begin{theorem} \label{campbellnullteo2}
Under hypothesis of Theorem \ref{campbellnullteo1}, if $X_t,Y_t$ are $\nabla^{G}$-martingales in $G$ with the null quadratic variation property, then $X_t\cdot Y_t$ is a $\nabla^{G}$-martingale in $G$.
\end{theorem}
\begin{proof}
Let $X_t,Y_t$ be $\nabla^{G}$-martingales in $G$. By Corollary \ref{ExpLogGroupcor2}, it is sufficient to show that $\mathcal{L}^{G\g}(X_t \cdot Y_t)$ is a  $\nabla^{\g}$-martingale.

From Theorem \ref{campbellnullteo1} we see that
\begin{equation}\label{martingalesproduct1}
\log{(X_{t} \cdot Y_{t})} = \int_0^t Ad(Y_s^{-1}) d\mathcal{L}^{G\mathfrak{g}}(X_s) + \log{(Y_t)}.
\end{equation}
Since $X_t,Y_t$ are $\nabla^{\g}$-martingales, Corollary \ref{ExpLogGroupcor2} assures that $\log{(X_t)}$ and $\log{(Y_t)}$ are $\nabla^{\g}$-martingales. Also, $\int_0^t Ad(Y_s^{-1}) d\mathcal{L}^{G\mathfrak{g}}(X_s)$ is a $\nabla^{\g}$-martingale because for any 1-form $\theta$ on $\g$ we have
\[
\!\!\int_0^t \!\!\theta d^{g}\!\!\int_0^s\!\!\! Ad(Y_r^{-1}) d^{g}\log{(X_r)}\!\! = \!\!\! \int_0^t \!\!\theta
Ad(Y_s^{-1}) d^{\g}\log{(X_s)}\!\! = \!\!\! \int_0^t \!\!\! (Ad(Y_s^{-1})^{*}\theta) d^{g}\log{(X_s)},
\]
Thus the sum of right side of (\ref{martingalesproduct1}) yields a $\nabla^{g}$-martingale and, consequently, the proof is complete. \qed
\end{proof}


In follows, we generalize the result due to P. Catuogno and P. Ruffino  \cite{catuogno1} for product of harmonic maps. Before, we are going to introduce a stochastic characterization for harmonic maps. Let $(M,g)$ be a Riemannian manifold and $N$ a smooth manifold with a connection $\nabla^{N}$. A smooth map $F:M \rightarrow N$ is a harmonic map if and only if it sends $g$-Brownian motions to $\nabla^{N}$-martingales.

\begin{proposition}\label{campbellnullteo3}
Let $(M_{j},g_{j})$, $j = 1,\ldots n$, be Riemannian manifolds, $G$ a Lie group with a complete, left-invariant connection $\nabla^{G}$ such that $\alpha(A,A) = 0$ for all $A \in \g$  and $\g$ its Lie algebra endowed with a complete connection $\nabla^{\g}$. If \linebreak $\phi_{j}:(M_{j},g_{j}) \rightarrow
(G,\nabla^{G})$ are  harmonic maps, then the product map $\phi_{1}\cdot\phi_{2}\cdot \ldots \cdot \phi_{n}$  from $M_{1}\times M_{2} \times \ldots \times M_{n}$ into $G$ is a harmonic map.
\end{proposition}
\begin{proof}
It is enough to take $n=2$. Let $\phi_{1}: (M_{1},g_{1}) \rightarrow (G, \nabla^{G})$ and \linebreak $\phi_{2}:(M_{2},g_{2}) \rightarrow (G, \nabla^{G})$ be harmonic maps.  Let $B^1_t$ and $B^2_t$ two independent Brownian motions in $M_{1}$ and $M_{2}$, respectively. Thus $(B^{1}_t,B^{2}_t)$ is a Brownian motion in the Riemmanian product manifold $M_{1} \times M_{2}$. It is sufficient to show that $\phi_1(B^{1}_t)\cdot\phi_2(B^{2}_t)$ is a $\nabla^{G}$-martingale. Since $B^{1}_t$ and $B^{2}_t$ are independent, $\phi_1(B^{1}_t)$ and $\phi_2(B^{2}_t)$ are too. Consequently, they have the null quadratic variation. Being $\phi_1(B^{1}_t)$ and $\phi_2(B^{2}_t)$ $\nabla^{G}$-martingales, Theorem \ref{campbellnullteo2} now shows that $\phi_1(B^{1}_t)\cdot\phi_2(B^{2}_t)$ is a $\nabla^{G}$-martingale, and the proof is complete. \qed
\end{proof}

\begin{example}
Let $G$ be a Lie group with a bi-invariant metric and $\g$ its Lie algebra with the Levi-Civita connection $\nabla^{G}$. Seeing the product on $G$ as application $m:(G\times G,\nabla^{G}\times \nabla^{G}) \rightarrow (G, \nabla^{G})$, Theorem \ref{campbellnullteo2} shows that the product $m$ is harmonic map.\qed
\end{example}

\begin{example}
Let $G$ be a Lie group equipped with a left-invariant connection $\nabla^{G}$ such that its connection function $\alpha$ satisfies $\alpha(A,A) = 0$ for all $A \in \g$ and $\g$ its Lie algebra endowed with a complete connection $\nabla^{g}$. Let $\gamma_{i}$ be $\nabla^{G}$-geodesics in $G$, $i=1, \ldots, n$. A map $f: (\mathbb{R}^{n},<,>) \rightarrow G$ defined by
\[
f(t_{1}, t_{2}, \ldots , t_{n}) = \gamma_{1}(t_{1}) \cdot \gamma_{1}(t_{2}) \cdot \ldots \cdot \gamma_{n}(t_{n}).
\]
is harmonic. It is a direct consequence of Theorem \ref{campbellnullteo3}. Indeed, it is sufficient to see any geodesic $\gamma_{i}$, $i=1, \ldots , n$, is a harmonic map.

In particular, assume that $G$ has a bi-invariant metric. Choose $n$ vectors $X_{1}, X_{2}, \ldots , X_{n} \in G$ such that $exp(t_{1}X_{1}), exp(t_{2}X_{2}), \ldots, exp(t_{n}X_{n})$ are geodesics (see for instant \cite{urakawa}). According to the facts above, $\exp(t_{1}X_{1}) \cdot \exp(t_{2}X_{2}) \cdot \ldots \cdot (\exp{t_{n}X_{n}})$ is a harmonic map. This example is also founded in \cite{urakawa} and \cite{catuogno1}.
\end{example}

\section{It\^o logarithm and examples}

We begin recalling the definition of the stochastic logarithm  due to M. Hakim-Dowek  and D. L\'epingle \cite{hakim}. Let $X_t$ be a semimartingale in $G$, then the stochastic logarithm, denoted by $L(X_{t})$, is the solution of the stochastic differential equation
\[
\delta L(X_{t})= L_{(X_{t})^{-1}*}(X_{t}) \delta X_{t} \ \ and \ \ X_{0}=e.
\]
Denoting by $\omega$ the Maurrer-Cartan form on $G$, it is simple to see that $L(X_{t}) = \int_{0}^{t}\omega \delta X_{s}$. Similarly, taking in account that the It\^o logarithm is a solution to stochastic differential equation (\ref{logaritmo}) we get $\log{(X_{t})} = \int_{0}^{t} \omega d^{\nabla^{G}}X_{s}$.

Following, we give a relation between the stochastic logarithm and It\^o logarithm.

\begin{proposition}\label{logmartingale}
Let $G$ be a Lie group with a complete, left invariant connection $\nabla^{G}$ with the associated connection function $\alpha$ and $\g$ its Lie algebra endowed with a complete connection $\nabla^{\g}$. A semimartingale $X_{t}$ in $G$ is a $\nabla^{G}$-martingale if and only if
\[
L(X_{t})+ \frac{1}{2}\int_{0}^{t}\alpha(L(X_{s}),L(X_{s})).
\]
is a local martingale in $\g$.
\end{proposition}
\begin{proof}
Let $X_{t}$ be a semimartingale in $G$. From the  It\^o-Stratonovich formula of conversion (\ref{conversion}) we compute
\begin{eqnarray*}
\log{(X_{t})}
& = & \int_{0}^{t} \omega_{G}d^{\nabla^{G}}X_{s} = \int_{0}^{t} \omega_{G}\delta X_{s} + \frac{1}{2}\int_{0}^{t}\nabla^{G}\omega(dX_{s},dX_{s})\\
& = & \int_{0}^{t} \omega_{G}\delta X_{s} + \frac{1}{2}\int_{0}^{t}\alpha(\omega dX_{s},\omega dX_{s})\\
& = & L(X_{t})+ \frac{1}{2}\int_{0}^{t}\alpha(L(X_{s}),L(X_{s})).
\end{eqnarray*}
Therefore the proof follows from Corollary \ref{ExpLogGroupcor2}.\qed
\end{proof}

The principal significance of this Proposition is that it allows us to see the geometry related whit martingales in Lie groups. Specifically, $\frac{1}{2}\int_{0}^{t}\alpha(L(X_{s}),L(X_{s}))$ is the term that differentiate the martingales in Lie Groups in accord to geometry given by the connection $\nabla^{G}$. We see this fact in the next two example.

\begin{example}
If $G$ has a bi-invariant metric, then the Levi-Civita connection on $G$ is given by $(\nabla^{G}_{A}B)(e) = \frac{1}{2}[A,B]$, where $A, B \in \g$. So the connection function $\alpha$ associate to $\nabla^{G}$ is $\alpha(A,B) = \frac{1}{2}[A,B]$. Therefore, $\alpha(A,A) = 0$. We thus conclude that $X_{t}$ is a $\nabla^{G}$-martingale if and only if $L(X_t)$ is a local martingale in $\g$.

A direct application of this is any semisimple Lie group $G$ equipped with the metric given by the Killing form.\qed
\end{example}

\begin{example}\label{leftmetricmartingale}
Suppose that $G$ is equipped with a complete left-invariant metric. Then the Levi-Civita connection on $G$ is given by
\[
(\nabla^{G}_{A}B)(e) = \frac{1}{2}[A,B] + U(A,B),
\]
where $A, B \in \g$ and $U: \g\times \g \rightarrow \g$ is a bilinear mapping defined by
\[
2<U(A,B),C> = <A,[C,B]> + <[C,A],B>,\ \ for\, all\, A,B,C\, \in \, \g.
\]
Here, $ <,> $ is the scalar product on $\g$ associated to the left metric on $G$. It follows that the connection function $\alpha$ associated to $\nabla^{G}$ is given by $\alpha(A,B) = \frac{1}{2}[A,B] + U(A,B)$. Consequently, $\alpha(A,A) = U(A,A)$. We thus conclude that $X_{t}$ is a $\nabla^{G}$-martingale if and only if $L(X_t)+ \frac{1}{2}\int_{0}^{t}U(L(X_{s}),L(X_{s}))$ is a local martingale in $\g$.\qed
\end{example}

In the sequel, we study the martingales in some matrix Lie groups with a left invariant metric. The idea is based in the work \cite{urakawa}, where it describes the bilinear mapping $U$ for some specific Lie groups. This fact and example above allows to characterize the martingales in these Lie groups. We begin with the Euclidian motion group $SE(3)$ and, in the sequel, with the three-dimensional non-compact Lie groups $SE(2), E(1,1), N^{3}$ and $SL(2,\mathbb{R})$.

\begin{example}[Euclidian motion group $SE(3)$]
The Euclidian motion group $SE(3)$ is defined by
\[
SE(3) =
\left\{
(X,u) =
\left(
\begin{array}{cc}
X & \mathbf{u}\\
\mathbf{0} & 0
\end{array}
\right); \mathbf{X} \in SO(3),\, \mathbf{u} \in \mathbb{R}^{n} \right\},
\]
where $\mathbf{0}$ is the $1\times 3$ matrix consisting of 0 and $\mathbf{u}$ is a 3-column vector in $\mathbb{R}^{3}$. The Lie algebra $\mathfrak{se}(3)$  is given by
\[
\mathfrak{se}(3) =
\left\{
(X,u) =
\left(
\begin{array}{cc}
X & \mathbf{u}\\
\mathbf{0} & 0
\end{array}
\right); \mathbf{X} \in \mathfrak{so}(3),\, \mathbf{u} \in \mathbb{R}^{n} \right\}.
\]
For our study we consider the inner product $<,>_{\lambda}$, with $\lambda>0$, defined by \linebreak $<(A,x),(B,y)>_{\lambda} = -\frac{1}{2}\mathrm{tr}(AB)+ \lambda^{2}x^ty$, where $x^{t}$ is the transpose of $x$. Take the basis $\beta=\{E_{1},E_{2},E_{3},e_{1},e_{2},e_{3}\}$ of $\mathfrak{se(3)}$, where

\begin{equation*}
E_{1}=\left(
\begin{array}{ccc}
0 & 0 & 0\\
0 & 0 & -1\\
0 & 1 & 0\\
\end{array}
\right),\,
E_{2}=\left(
\begin{array}{ccc}
0 & 0 & 1 \\
0 & 0 & 0\\
-1 & 0 & 0\\
\end{array}
\right),\,
E_{3}=\left(
\begin{array}{ccc}
0 & -1 & 0\\
1 & 0 & 0\\
0 & 0 & 0\\
\end{array}
\right)
\end{equation*}
and
\begin{equation*}
e_{1}=\left(
\begin{array}{c}
1 \\
0 \\
0 \\
\end{array}
\right),\,
e_{2}=\left(
\begin{array}{ccc}
0 \\
1 \\
0 \\
\end{array}
\right),\,
e_{3}=\left(
\begin{array}{ccc}
0\\
0\\
1\\
\end{array}
\right).
\end{equation*}

Let $X_{t}$ be a semimartingale in $SE(3)$ and $L(X_{t})$ the stochastic logarithm in $\mathfrak{se}(3)$. According to basis $\beta$, we may write $L(X_{t})= \sum_{1}^{3} x_{i}(t)E_{i} + \sum_{1}^{3} y_{i}(t)e_{i}$, where $x_{i}(t), y_{i}(t)$ are real semimartingales. Then from Lemma 5.1 in \cite{urakawa} we see that $U(L(X_t),L(X_{t})) = x(t)\times y(t)$(the vector product), where
\begin{equation*}
x(t)=\sum_{i=1}^{3}x_{i}(t) =
\left(
\begin{array}{c}
x_{1}(t)\\
x_{2}(t)\\
x_{3}(t)\\
\end{array}
\right) \in \mathbb{R}^{3} \, and \,
y(t)=\sum_{i=1}^{3}y_{i}(t) =
\left(
\begin{array}{c}
y_{1}(t)\\
y_{2}(t)\\
y_{3}(t)\\
\end{array}
\right) \in \mathbb{R}^{3}.
\end{equation*}
From Example \ref{leftmetricmartingale} we see that $X_{t}$ is a $\nabla^{G}$-martingale if and only if
\[
L(X_{t}) + \frac{1}{2}\int_{0}^{t} x(s)\times y(s)
\]
is a local martingale in $\mathfrak{se}(3)$.\qed
\end{example}

\begin{example}[SE(2)]
With a little modification of the example above, we have as a basis for the Lie algebra $\mathfrak{se}(2)$
\begin{equation*}
H=\left(\left(
\begin{array}{cc}
0 & -1\\
1 & 0
\end{array}
\right)
,
\left(
\begin{array}{c}
0 \\
0
\end{array}
\right)\right),
\,
e_{1}=\left(
\mathbf{0}, \,
\left(
\begin{array}{c}
1 \\
0
\end{array}
\right)\right),\,
e_{2}=\left(
\mathbf{0}, \,
\left(
\begin{array}{c}
0 \\
1
\end{array}
\right)\right).
\end{equation*}
The inner product adopt $<,>_{\lambda}$, with $\lambda>0$, in $\mathfrak{se}(2)$ is defined by
\[
<(aH + x_{1}e_{1}+x_{2}e_{2}),(bH+y_{1}e_{1}+y_{2}e_{2})>  = ab +\lambda^{2}(x_{1}y_{1}+x_{2}y_{2}).
\]
Let $X_{t}$ be a semimartingale in $SE(2)$ and $L(X_{t})$ the stochastic logarithm in $\mathfrak{se}(2)$. We may write $ L(X_{t})= a(t)H + a_{1}(t)e_{1}+a_{2}(t)e_{2}$, where $a(t),a_{1}(t), a_{2}(t)$ are real semimartingales. Lemma 6.1 in \cite{urakawa} now assures that $U(L(X_t),L(X_{t})) = a(t)H(a_{1}(t)e_{1}+ a_{2}(t)e_{2})$. We conclude from Example \ref{leftmetricmartingale} that $X_{t}$ is a $\nabla^{G}$-martingale if and only if
\[
L(X_{t}) + \frac{1}{2}\int_{0}^{t} a(s)H(a_{1}(s)e_{1}+ a_{2}(s)e_{2})
\]
is a local martingale in $\mathfrak{se}(2)$.\qed
\end{example}


\begin{example}[E(1,1)]
The three-dimensional Lie group $E(1,1)$ is given by
\[
E(1,1)=
\left\{
\left(
\begin{array}{ccc}
\exp(\xi) & 0 & x_1\\
0 & \exp(-\xi) & x_{2}\\
0 & 0 & 1
\end{array}
\right);
\, \xi,x_{1},x_{2} \in \mathbb{R}
\right\},
\]
with standard multiplication of matrix. Its Lie algebra $\mathfrak{e}(1,1)$ is given by
\[
\mathfrak{e}(1,1)=
\left\{
\left(
\left(
\begin{array}{cc}
\xi & 0\\
0 & -\xi
\end{array}
\right),\,
\left(
\begin{array}{c}
x_{1}\\
x_{2}
\end{array}
\right)
\right)
=
\left(
\begin{array}{ccc}
\xi & 0 & x_1\\
0 & -\xi & x_{2}\\
0 & 0 & 1
\end{array}
\right);
\, \xi,x_{1},x_{2} \in \mathbb{R}
\right\}.
\]
A basis for $\mathfrak{e}(1,1)$ is
\begin{equation*}
H=\left(\left(
\begin{array}{cc}
1 & 0\\
0 & -1
\end{array}
\right)
,
\left(
\begin{array}{c}
0 \\
0
\end{array}
\right)\right),
\,
e_{1}=\left(
\mathbf{0}, \,
\left(
\begin{array}{c}
1 \\
0
\end{array}
\right)\right),\,
e_{2}=\left(
\mathbf{0}, \,
\left(
\begin{array}{c}
0 \\
1
\end{array}
\right)\right).
\end{equation*}
The inner product $<,>_{\lambda}$ $(\lambda>0)$ adopted on $\mathfrak{e}(1,1)$ is
\[
<(aH + x_{1}e_{1}+x_{2}e_{2}),(bH+y_{1}e_{1}+y_{2}e_{2})>_{\lambda}  = ab +\lambda^{2}(x_{1}y_{1}+x_{2}y_{2}).
\]
Take a semimartingale $X_{t}$ in $E(1,1)$ and the stochastic logarithm $L(X_{t})$ in $\mathfrak{e}(1,1)$. So the stochastic logarithm may be written as  $ L(X_{t})= a(t)H + a_{1}(t)e_{1}+a_{2}(t)e_{2}$. Then from Lemma 6.4 in \cite{urakawa} we see that $U(L(X_t),L(X_{t})) = \|a_{1}(t)e_{1}+a_{2}(t)e_{2}\|^{2} \lambda^{2}H - a(t)H(a_{1}(t)e_{1}+ a_{2}(t)e_{2})$, where $\|a_{1}(t)e_{1}+a_{2}(t)e_{2}\|^{2} = a_{1}^{2}(t)-a_{2}^{2}(t)$. Therefore, Example \ref{leftmetricmartingale} assures that $X_{t}$ is a $\nabla^{G}$-martingale if and only if
\[
L(X_{t}) + \frac{1}{2}\int_{0}^{t} \|a_{1}(s)e_{1}+a_{2}(s)e_{2}\|^{2} \lambda^{2}H - a(s)H(a_{1}(s)e_{1}+ a_{2}(s)e_{2}))
\]
is a local martingale in $\mathfrak{e}(1,1)$.
\end{example}

\begin{example}[$N^{3}$]
The Heisenberg group is the three dimensional nilpotent group defined by
\[
N^{3}=
\left\{
\left(
\begin{array}{ccc}
1 & x & z\\
0 & 1 & y\\
0 & 0 & 1
\end{array}
\right);
\, x,y,z \in \mathbb{R}
\right\}.
\]
The Lie algebra $\mathfrak{n}^{3}$ of the Heisenberg group is
\[
\mathfrak{n}^{3}=
\left\{
\left(
\begin{array}{ccc}
0 & x & z\\
0 & 0 & y\\
0 & 0 & 0
\end{array}
\right);
\, x,y,z \in \mathbb{R}
\right\}.
\]
A basis for $\mathfrak{n}^{3}$ is
\[
X=\left(
\begin{array}{ccc}
0 & 1 & 0\\
0 & 0 & 0\\
0 & 0 & 0
\end{array}
\right),\,
Y=\left(
\begin{array}{ccc}
0 & 0 & 0\\
0 & 0 & 1\\
0 & 0 & 0
\end{array}
\right),\,
Z=\left(
\begin{array}{ccc}
0 & 0 & 1\\
0 & 0 & 0\\
0 & 0 & 0
\end{array}
\right).
\]
We adopt the inner product $<,>_{\lambda}$ $(\lambda>0)$  on $\mathfrak{n}^{3}$ given by
\[
<(a_1X + b_1Y+ c_1Z),(a_2X+b_2Y+c_2Z)>  = a_1a_2 +\lambda^{2}(b_{1}b_{2}+c_{1}c_{3}).
\]
Let $X_{t}$ be a semimartingale in $N^{3}$ and $L(X_{t})$ the stochastic logarithm in $\mathfrak{n}^{3}$. Write $ L(X_{t})= a(t)H + b(t)E_{+}+ c(t)E_{-}$ in terms of basis above. Using Lemma 6.9 in \cite{urakawa} we obtain $U(L(X_t),L(X_{t})) = \lambda^{2}b(t)c(t)X +\lambda^{2}a(t)c(t)Y$. So Example \ref{leftmetricmartingale} shows that $X_{t}$ is a $\nabla^{G}$-martingale if and only if
\[
L(X_{t}) + \frac{1}{2}\int_{0}^{t} -\lambda^{2}b(s)c(s)X +\lambda^{2}a(s)c(s)Y
\]
is a local martingale in $\mathfrak{n}^{3}$.
\end{example}

\begin{example}[$SL(2,\mathbb{R})$]
The Lie algebra $\mathfrak{sl}(2,\mathbb{R})$ of $SL(2,\mathbb{R})$ is $\{H,E_{+},E_{-}\}$, where
\[
H=
\left(
\begin{array}{cc}
1 & 0 \\
0 & -1
\end{array}
\right),\,
E_{+}=
\left(
\begin{array}{cc}
0 & 1 \\
0 & 0
\end{array}
\right),\,
E_{-}=
\left(
\begin{array}{cc}
0 & 0 \\
1 & 0
\end{array}
\right).
\]
Assume that $\mathfrak{sl}(2,\mathbb{R})$ is equipped with the inner product $<,>_{\lambda}$ $(\lambda>0)$  given by
\[
<(a_1H + b_1E_{+}+ c_1E_{-}),(a_2H+b_2E_{+}+c_2E_{-})>  = a_1a_2 +\lambda^{2}(b_{1}b_{2}+c_{1}c_{3}).
\]
Taking a semimartingale $X_{t}$  in $SL(2,\mathbb{R})$ and the stochastic logarithm $L(X_{t})$  in $\mathfrak{sl}(2,\mathbb{R})$ we have $ L(X_{t})= a(t)X + b(t)Y+ c(t)Z$, where $a(t), b(t), c(t)$ are real semimartingales. Applying Lemma 6.4 in \cite{urakawa} we see that $U(L(X_t),L(X_{t})) =  \frac{2}{\lambda^{2}}(b(t)^{2}-c(t)^{2})H + (-2a(t)b(t)+a(t)c(t)\lambda)E_{+} + (-a(t)b(t)\lambda +2a(t)c(t))E_{-}$. Now, Example \ref{leftmetricmartingale} assures that $X_{t}$ is a $\nabla^{G}$-martingale if and only if
\[
\begin{array}{l}
L(X_{t}) + \\
 \frac{1}{2}\displaystyle\int_{0}^{t}\!\!\! \displaystyle\frac{2}{\lambda^{2}}(b(s)^{2}-c(s)^{2})H\! +\! (-2a(s)b(s)+a(s)c(s)\lambda)E_{+}\! +\! (-a(s)b(s)\lambda +2a(s)c(s))E_{-}
\end{array}
\]
is a local martingale in $\mathfrak{sl}(2,\mathbb{R})$.

\end{example}

%

\end{document}